\documentclass[11pt,a4paper]{article}

\usepackage{latexsym}
\usepackage{amssymb}

\newtheorem{theorem}{Theorem}[section]
\newtheorem{lemma}[theorem]{Lemma}
\newtheorem{remark}[theorem]{Remark}
\newtheorem{example}[theorem]{Example}
\newtheorem{proposition}[theorem]{Proposition}
\newtheorem{definition}[theorem]{Definition}
\newtheorem{corollary}[theorem]{Corollary}
\newtheorem{assumption}[theorem]{Assumption}

\newcommand{\bbr}{{\mathbb R}}  
\newcommand{\Expect}[1]{{\rm E}\left[\,#1\,\right]}
\renewcommand{\P}[1]{{\rm P}\left[#1\right]}
\newcommand{\defeq}{\stackrel{\rm def}{=}}
\newcommand{\bs}{\backslash}
\newcommand{\ts}[1]{\textstyle{#1}}
\newcommand{\lint}{\int\limits}
\newcommand{\M}{{\cal M}} 
\newcommand{\BB}{{\cal B}}

\title{Conditional Expectation as Quantile Derivative}
\author{
Dirk Tasche\thanks{Zentrum Mathematik (SCA), TU M\"unchen, 80290 M\"unchen,
Germany;
email: tasche@ma.tum.de}
}

\date{November 13, 2000}

\begin{document}
\maketitle

\begin{abstract}
For a linear combination $\sum u_j\,X_j$ of random variables,
we are interested in the partial derivatives of its $\alpha$-quantile $Q_\alpha(u)$ 
regarded as a function of the weight vector $u = (u_j)$. It turns out that under 
suitable conditions on the joint distribution of $(X_j)$ the derivatives exist
and coincide with the conditional expectations of the $X_i$ given that $\sum u_j\,X_j$
takes the value $Q_\alpha(u)$. Moreover, using this result, we deduce formulas for
the derivatives with respect to the $u_i$ for the so-called expected shortfall
$\Expect{\big|\sum u_j\,X_j - Q_\alpha(u)\big|^\delta\,\big|\,\sum u_j\,X_j \le Q_\alpha(u)}$,
with $\delta \ge 1$ fixed. Finally, we study in some more detail the coherence properties
of the expected shortfall in case $\delta =1$.

{\sc Key words:} \textit{quantile; value-at-risk; quantile derivative; 
conditional expectation; expected shortfall; conditional value-at-risk; coherent
risk measure.}

\end{abstract}

\section{Introduction}
\label{sec:1}

The last decade has seen a growing interest in quantiles of
probability distributions by practitioners mainly in the financial
industry. Since quantiles have a  simple
interpretation in terms of over- or under\-shoot probabilities
they have found entrance in current risk management practice
in form of the value-at-risk concept (cf. \cite{Jorion}).

In particular, there is need for computing derivatives of 
quantiles of weighted sums of
random variables with respect to the weights. \cite{Arzac_Bawa}
represents an early example for the use of these derivatives.
More recently, in \cite{Tasche99} was shown for general
risk measures that their derivatives with respect to the 
asset weights are the key to the
solution of the ``capital allocation'' problem (see also \cite{Delbaen_Denault}
for the case of ``coherent'' risk measures). The problem to
allocate the total risk to risk sources in connection with
the need to differentiate risk measures appears also in other
scientific disciplines. For an example in statistics see~\cite{Tasche00}.

In case of normally distributed random vectors the formulae
for the derivatives are obvious (cf. \cite{Garman97}).
In \cite{Mal_Kib_Che}, a result (Theorem 1) was provided  for general
distributions and even non-linear combinations of random
variables. A similar formula can be found in \cite{Gaiv_Pflug}.
Unfortunately, the  results in \cite{Mal_Kib_Che} and \cite{Gaiv_Pflug}
lack of an intuitive explanation.
Nevertheless, readily interpretable formulae are available
for the case of linear combinations of random variables
(cf. \cite{Hallerbach99} or \cite{Scaillet}).
The primary intention with the present paper is to 
give
a sufficient condition as general as possible on the underlying distribution 
for the formulae in  \cite{Hallerbach99} and \cite{Scaillet}) to remain
valid. In addition, we will transfer the result on quantiles onto the so-called
expected shortfall (also called conditional value-at-risk). Finally,
because of the theoretical importance of the expected shortfall we will
discuss its role as coherent risk measure in the sense of \cite{ADEH99}.

This paper is organized as follows: In section \ref{sec:2} we recall
some properties of conditional densities and introduce the
technical assumptions needed for the main result on differentiation
of quantiles. In section \ref{sec:3}, this result  is presented in
an easy to digest (eq. (\ref{eq:3.4})) and a rigorous version
(Theorem \ref{th:2}). We then apply the result to the expected shortfall
in section \ref{sec:4}. The last section is devoted to a
more detailed study of the latter with respect to its coherence 
properties. 

\section{Some facts on conditional densities}
\label{sec:2}

We are going to present sufficient conditions for quantiles
of a sum $\sum u_i\,X_i$ to be differentiable with respect to
the weights $u_i$. These conditions will heavily rely on the
existence of a conditional density of one component of the
$X_i$ given the others. So we will start our study by summarizing
some facts on conditional densities. 

First we recall the notion
of conditional density in the context of a random vector
$(X_1, \ldots, X_d)$. We write $\BB(\bbr^d)$ for the $\sigma$-algebra of
Borel sets on $\bbr^d$. By the indicator function $I(A) = I(A, \omega)$ of a set $A$ we
mean the function defined by
\begin{equation}
I(A,\omega)\quad  \defeq\quad \left\{\begin{array}{rl} 
1, & \mbox{if\ } \omega \in A\\
0, & \mbox{if\ } \omega \notin A\,.
\end{array} \right.
  \label{eq:2.0}
\end{equation}
\begin{definition}
  \label{de:2.1}
Let $d \ge 2$ and let $(X_1, \ldots, X_d)$ be an $\bbr^d$-valued random vector.
A measurable function $\phi: \bbr^d \to [0,\infty), (t, x_2, \ldots, x_d) \mapsto \phi(t, x_2, \ldots,x_d)$
is called {\em conditional density} of $X_1$ given $(X_2, \ldots, X_d)$ if for all $A \in \BB(\bbr),\
B \in \BB(\bbr^{d-1})$
$$
\P{X_1 \in A,\ (X_2, \ldots, X_d)\in B}\quad = \quad 
\Expect{ I((X_2, \ldots, X_d)\in B) \int_A \phi(t, X_2, \ldots, X_d) dt }.
$$
\end{definition}

An equivalent formulation for Definition \ref{de:2.1} is
\begin{equation}
 \int_A \phi(t, X_2, \ldots, X_d)\, dt 
 \quad =\quad \P{ X_1 \in A \,\Big|\,X_2, \ldots, X_d},
  \label{eq:2.1}
\end{equation}
i.e. $\phi(\cdot, X_2,\ldots, X_d)$ is a density of the conditional
distribution of $X_1$ given $(X_2, \ldots, X_d)$. 
Recall the well-known fact that the existence of a joint
density of $(X_1, \ldots, X_d)$ is sufficient but not necessary for
the existence of a conditional density of $X_1$ given $(X_2, \ldots, X_d)$.
On the other hand, the existence of such a conditional density implies
that the unconditional distribution of $X_1$ has a density $f$ that is given
by $f(t) = \Expect{\phi(t, X_2, \ldots, X_d)}$. 
Moreover, a situation can occur where the distribution of $(X_2, \ldots, X_d)$ is
purely discrete and a conditional density of $X_1$ given $(X_2, \ldots, X_d)$ exists.

For our purpose, the following three easy conclusions from the existence of a conditional
density are important.
\begin{lemma}
  \label{le:2.2} 
Let $d \ge 2$ and let $(X_1, \ldots, X_d)$ be an $\bbr^d$-valued random vector
with a conditional density $\phi$ of $X_1$ given $(X_2, \ldots, X_d)$. Then for any 
weight vector $(u_1, \ldots, u_d) \in \bbr\bs\{0\}\times \bbr^{d-1}$ we have
\begin{enumerate}
\item the function $t \mapsto |u_1|^{-1}\, \Expect{ \phi\left( u_1^{-1}\left( t -\ts{ \sum_{j=2}^d u_j\,X_j}\right),
X_2, \ldots, X_d\right) }$ is a density of\\
$\sum_{i=1}^d u_i\,X_i\,$,
\item for $i = 2, \ldots, d$, almost surely for $t\in\bbr$
$$
\Expect{X_i\,\bigg|\,\sum_{j=1}^d u_j\,X_j = t} =
\frac{\Expect{ X_i\, \phi\left( u_1^{-1}\left( t - \ts{\sum_{j=2}^d u_j\,X_j}\right),
X_2, \ldots, X_d\right) }}
{\Expect{ \phi\left( u_1^{-1}\left( t -\ts{ \sum_{j=2}^d u_j\,X_j}\right),
X_2, \ldots, X_d\right) }}\,,\quad\mbox{and}
$$
\item almost surely for $t\in\bbr$
$$
\Expect{X_1\,\bigg|\,\sum_{j=1}^d u_j\,X_j = t} = \frac
{\Expect{ \frac{t-\sum_{i=2}^d u_i\,X_i}{u_1}\, \phi\left( u_1^{-1}\left( t - \ts{\sum_{j=2}^d u_j\,X_j}\right),
X_2, \ldots, X_d\right) }}
{\Expect{ \phi\left( u_1^{-1}\left( t -\ts{ \sum_{j=2}^d u_j\,X_j}\right),
X_2, \ldots, X_d\right) }}\,.
$$ 
\end{enumerate}
\end{lemma}
To say Lemma \ref{le:2.2} with words: if there is a conditional density of $X_1$ given the other
components, then subject to the condition $u_1 \not= 0$ the distribution of $\sum_{i=1}^d u_i\, X_i$ 
is absolutely continuous with 
density specified in Lemma \ref{le:2.2} (i), and the conditional expectations of the $X_i$ given the
sum $\sum_{j=1}^d u_j\, X_j$ can be 
calculated via the formulae in Lemma \ref{le:2.2} (ii) and (iii).

In the subsequent section, it will turn out that the quantiles of the sum $\sum_{i=1}^d u_i\, X_i$
are differentiable with respect to the weights $u_i$ if the quantities mentioned in Lemma \ref{le:2.2}
are smooth in a certain sense. This observation motivates the following definition.

\begin{assumption}
  \label{as:2.3}
Let $d \ge 2$ and let $(X_1, \ldots, X_d)$ be an $\bbr^d$-valued random vector
with a conditional density $\phi$ of $X_1$ given $(X_2, \ldots, X_d)$.
We say that $\phi$
satisfies  Assumption \ref{as:2.3} in an open set $U \subset \bbr\bs\{0\}\times \bbr^{d-1}$ if
the  following three conditions hold:
\begin{enumerate}
\item 
For fixed $x_2, \ldots, x_d$ the function $t \mapsto \phi(t, x_2, \ldots, x_d)$
is continuous in $t$.
\item
The mapping 
$$
(t, u) \mapsto \Expect{ \phi\left( u_1^{-1}\left( t -\ts{ \sum_{j=2}^d u_j\,X_j}\right),
X_2, \ldots, X_d\right) }, \quad 
\bbr \times U \to [0,\infty)
$$
is finite-valued and continuous.
\item 
For each $i = 2, \ldots, d$ the mapping
$$
(t, u) \mapsto \Expect{ X_i\, \phi\left( u_1^{-1}\left( t - \ts{\sum_{j=2}^d u_j\,X_j}\right),
X_2, \ldots, X_d\right) }, \quad 
\bbr \times U \to \bbr
$$
is finite-valued and continuous.
\end{enumerate}
\end{assumption}

Note that (i) from Assumption \ref{as:2.3} in general does imply neither (ii) nor (iii). Furthermore,
(ii) and (iii) may be valid even if the components of the random
vector $(X_1,\ldots, X_d)$ do not have finite expectations. 

\begin{remark}
Here is a list of some situations in which Assumption \ref{as:2.3} is
satisfied:
\begin{itemize}
\item[1)] 
$(X_1,\ldots,X_d)$ is normally distributed and its covariance matrix has full rank.
\item[2)]
$(X_1,\ldots,X_d)$ and $\phi$ satisfy (i)
and for each $(s, v)\in 
\bbr\times U$ there is some neighbourhood $V$
such that the 
random fields 
$$
\Big( \phi\Big( u_1^{-1}( t - \textstyle{\sum_{j=2}^d u_j\,X_j}),
X_2, \ldots, X_d\Big) \Big)_{(t,u)\in V}
$$
and for $i= 2, \ldots, d$
$$
\Big(X_i\, \phi\Big( u_1^{-1}( t -\textstyle{ \sum_{j=2}^d u_j\,X_j}),
X_2, \ldots, X_d\Big) \Big)_{(t,u)\in V}
$$
are uniformly integrable.
\item[3)]
$\Expect{ \mid X_i\mid} < \infty, \ i=2,\ldots, d$, and 
$\phi$ is bounded and  satisfies (i).
\item[4)]
$\Expect{ \mid X_i\mid} < \infty, \ i=2,\ldots, d$.
$X_1$ and $(X_2,\ldots, X_d)$ are independent. $X_1$ has a 
continuous density.
\item[5)]
There is a finite set $M \subset \bbr^{d-1}$ such 
that
$\P{(X_2, \ldots, X_d) \in M } = 1\,,$
and (i) is satisfied.
\end{itemize}
  \label{rm:2.4}
\end{remark}
Note that Remark \ref{rm:2.4} 3)  is a special case of 2) and that 4) and 5) resp.~are
special cases of 3). Perhaps, 4) is the case most interesting for
applications. It corresponds to the situation where a sample of
$(X_1,\ldots,X_d)$ and a weight vector $u$ are given and the density
of $\sum_{j=1}^d u_j\,X_j$ is estimated by kernel estimation.

\section{Quantile Derivatives}
\label{sec:3}

If $X$ is a real valued random variable and $\alpha$ is any number between 0 and
1, the $\alpha$-quantile of $X$ is the $100 \alpha$\%-threshold of $X$, i.e. 
the lowest bound to be exceeded by $X$ only with probability $100 (1-\alpha)$\%.
We will make use of the following formal definition.

\begin{definition}   \label{de:1}
Let $X$ be a real valued random variable and let $\alpha \in (0,1)$.
Then the $\alpha${\em -quantile} $Q_\alpha(X)$ of $X$ is defined by
$$
Q_\alpha(X)\quad \defeq \quad \inf\{ x\in \bbr\,|\,\P{X \le x} \ge \alpha \}\,.
$$
\end{definition}
In general, the case $\P{X \le Q_\alpha(X)} > \alpha$ is possible, 
but in this paper solely  $\P{X \le Q_\alpha(X)} = \alpha$ will occur.
The reason is that our method for proving differentiability of the quantiles will be
based on the implicit function theorem.

Let us briefly outline the reasoning. We want to study the mapping
\begin{equation}
Q_\alpha(u) \quad \defeq \quad Q_\alpha\Big( \sum_{j=1}^d u_j\,X_j \Big)\,,
  \label{eq:3.1}
\end{equation}
regarded as a function of the weight vector $u$. Assume for the moment that
we already know that $Q_\alpha(u)$ is differentiable with respect to the
components of $u$. If there is a conditional density of $X_1$ given
$(X_2, \ldots, X_d)$, then by Lemma \ref{le:2.2} (i) the distribution of
$\sum_{j=1}^d u_j\,X_j$ is continuous, and we obtain by (\ref{eq:2.1}) for all $u$ with $u_1 > 0$
\begin{equation}
\alpha\  =\  \ts{\P{ \sum_{j=1}^d u_j\,X_j \le Q_\alpha(u)}}
\ =\  \mbox{E}\bigg[\lint_{-\infty}^{u_1^{-1} (Q_\alpha(u) - \sum_{j=2}^d u_j\,X_j)} 
\phi(t, X_2, \ldots, X_d)\, d t\bigg]\,.
  \label{eq:3.2}
\end{equation}
Ignoring the question whether or not differentiation under the expectation
is permitted, by differentiating with respect to $u_i, i= 2, \ldots, d,$ we obtain from
(\ref{eq:3.2})
\begin{equation}
0 \quad = \quad u_1^{-1} \, \Expect{ \Big( \frac{\partial Q_\alpha(u)}{\partial u_i}-X_i\Big)\, \phi\Big( u_1^{-1} 
(Q_\alpha(u) - \ts{\sum_{j=2}^d u_j\,X_j}), X_2, \ldots, X_d\Big)}\,.
  \label{eq:3.3}
\end{equation}
Solving (\ref{eq:3.3}) for $\frac{\partial Q_\alpha(u)}{\partial u_i}$ and applying 
formally Lemma 
\ref{le:2.2} (ii) now yields
\begin{equation}
\frac{\partial Q_\alpha(u)}{\partial u_i} \ =\ 
\Expect{ X_i\,\Big|\,\ts{\sum_{j=1}^d u_j\,X_j} = Q_\alpha(u)}\,,\ i= 2, \ldots, d. 
  \label{eq:3.4}
\end{equation}
An analogous computation could be done in the cases $u_1 < 0$ and $i=1$ and would yield 
(\ref{eq:3.4}) also for $u_1 < 0$ or $i=1$.
Equation (\ref{eq:3.4}) has been presented in \cite{Hallerbach99} without
examination of the question whether $Q_\alpha$ is differentiable
and in \cite{Scaillet} for the case of $(X_1, \ldots, X_d)$ with a
joint density.

In order to make this approach mathematically rigorous by invoking the
implicit function theorem, we have to verify some smoothness conditions
for the expression $\P{ \sum_{j=1}^d u_j\,X_j \le z}$ considered as
a function of $u$ and $z$.

\begin{lemma}
Let $(X_1,\ldots, X_d)$ be an $\bbr^d$-valued random vector.
Assume that there is a conditional density $\phi$ of 
$X_1$ given $(X_2, \ldots, X_d)$ satisfying Assumption \ref{as:2.3}
in some open set $U \subset
\bbr\bs\{0\} \times \bbr^{d-1}$.
Define the random field $(Z(u))_{u\in U}$ by 
\begin{equation}
Z(u) \quad \defeq \quad \sum_{i=1}^d u_i\,X_i\,, \quad u \in U\,.
\label{eq:2}
\end{equation}
Then the function $F: \bbr \times U \to [0, \infty)$, defined by
$$
  F(z, u) \quad\defeq \quad \P{ Z(u) \le z }\,,
$$
is partially differentiable in $z$ and $u_i$, $i = 1,\ldots, d$,
with jointly continuous derivatives
\begin{eqnarray}
\frac{\partial F(z,u)}{\partial z}  & = & 
|u_1|^{-1}\, \Expect{ \phi\bigg( u_1^{-1}\bigg( z - \sum_{j=2}^d u_j\,X_j\bigg),
X_2, \ldots, X_d\bigg) },\label{eq:4}\\[1ex]
\frac{\partial F(z,u)}{\partial u_1}  & = & 
\frac{-{\rm sign}(u_1)}{u_1^2}
\, \Expect{ \bigg(z- \sum_{j=2}^d u_j\,X_j\bigg) 
\phi\bigg( u_1^{-1}\bigg( z - \sum_{j=2}^d u_j\,X_j\bigg),
X_2, \ldots, X_d\bigg) },\label{eq:5}\\
 \mbox{and} & & \nonumber\\
\frac{\partial F(z,u)}{\partial u_i}  & = & 
-\,{|u_1|^{-1}}
\, \Expect{ X_i\,
\phi\bigg( u_1^{-1}\bigg( z - \sum_{j=2}^d u_j\,X_j\bigg),
X_2, \ldots, X_d\bigg) },\ i = 2, \ldots, d.\label{eq:6}
\end{eqnarray}
  \label{le:1}
\end{lemma}
\textbf{Proof.}
The joint continuity of the expressions for the
partial derivatives follows from Assumption \ref{as:2.3} (ii) and (iii).
Equation (\ref{eq:4}) is obvious since by Lemma \ref{le:2.2} (i) 
the right-hand side of 
(\ref{eq:4}) as function of $z$ is a continuous density of $Z(u)$.

By the representations
$$
  F(z, u) \quad =\quad {\rm E}\bigg[\,\lint_{-\infty}^{u_1^{-1}(z-\sum_{j=2}^d u_j\,X_j) }
\phi(t, X_2,\ldots, X_d)\,dt\,\bigg]
$$
in case $u_1 > 0$ and
$$
  F(z, u) \quad =\quad {\rm E}\bigg[\, \lint^{\infty}_{
u_1^{-1}(z-\sum_{j=2}^d u_j\,X_j)} 
\phi(t, X_2,\ldots, X_d)\,dt\,\bigg]
$$
in case $u_1 < 0$ respectively, the application
of  Theorem A.(9.1) from  \cite{Durrett} on  differentiation under the integral
yields the desired formulae (\ref{eq:5}) and (\ref{eq:6}). \hfill $\Box$

With Lemma \ref{le:1} we are in a position suitable to give a rigorous
formulation to (\ref{eq:3.4}). Keep in mind that by Lemma \ref{le:2.2}
equation
(\ref{eq:3.4}) on the one hand and (\ref{g211}) and (\ref{g212}) on the other
hand have
essentially the same meaning.

\begin{theorem}
Let $\alpha \in (0,1)$ be fixed,
and let $(X_1,\ldots, X_d)$ be an $\bbr^d$-valued random vector 
with a conditional density $\phi$ of 
$X_1$ given $(X_2, \ldots, X_d)$
that satisfies 
Assumption \ref{as:2.3} in some open set $U \subset 
\bbr\bs\{0\} \times \bbr^{d-1}$.
Define the random field $(Z(u))_{u\in U}$ by (\ref{eq:2}) and 
the function 
$
Q_\alpha: U \to \bbr
$ 
by (\ref{eq:3.1}).

If the density 
$t \mapsto |u_1|^{-1}\, \Expect{ \phi\left( u_1^{-1}\left( t -\ts{ \sum_{j=2}^d u_j\,X_j}\right),
X_2, \ldots, X_d\right) }$ of $Z(u)$ is positive at $t = Q_\alpha(u)$, then $Q_\alpha$
is partially differentiable in some neighbourhood of $u$ 
with continuous derivatives
\begin{eqnarray}
\lefteqn{\hspace{1cm}
\frac{\partial Q_\alpha}{\partial u_1}(u)  = }\label{g211} \\
& & 
u_1^{-1}\left(Q_\alpha(u) - 
\frac{\ts{\Expect{\left(\sum_{j=2}^d u_j\,X_j\right) 
\phi\left( u_1^{-1}(Q_\alpha(u) -\sum_{j=2}^d u_j\,X_j), X_2, \ldots, X_d\right) }}}
{\ts{\Expect{
\phi\left( u_1^{-1}(Q_\alpha(u) -\sum_{j=2}^d u_j\,X_j), X_2, \ldots, X_d\right) }}}
\right) \nonumber
\end{eqnarray}
and 
\begin{equation}\label{g212}\quad
\frac{\partial Q_\alpha}{\partial u_i}(u)  =  
\frac{\ts{\Expect{X_i\,
\phi\left( u_1^{-1}(Q_\alpha(u) -\sum_{j=2}^d u_j\,X_j), X_2, \ldots, X_d\right) }}}
{\ts{\Expect{
\phi\left( u_1^{-1}(Q_\alpha(u) -\sum_{j=2}^d u_j\,X_j), X_2, \ldots, X_d\right) }}},
\ i=2,\ldots,d.
\end{equation}
  \label{th:2}
\end{theorem}
\textbf{Proof.}
As Assumption \ref{as:2.3} implies that the distributions of
$(Z(u))_{u\in U}$ are continuous, we have for all $u \in U$
\begin{equation}
  \label{eq:9}
  \P{Z(u) \le Q_\alpha(u)}\quad =\quad \alpha\,.
\end{equation}
From (\ref{eq:9}), by Lemma \ref{le:1} and the implicit function theorem we
obtain the assertion. \hfill $\Box$

\section{Shortfall Derivatives}
\label{sec:4}

As a quantile at a fixed level gives only local information about the
underlying distribution, a promising way to escape from this
shortcoming is to consider the so-called expected shortfall 
over or under the quantile. These quantities, to be defined in the
following theorem, can be interpreted as moments of the difference between
the underlying
random sum and a quantile in a worst case situation specified by the 
confidence level of the quantile.

\begin{theorem}
Let $\delta \ge 1$ be fixed. Let $(X_1, \ldots, X_d)$, $\phi$, $U$, and $\alpha$ be as in Theorem \ref{th:2},
and assume 
additionally
$$
\Expect{\mid X_i\mid^\delta} < \infty,\quad i = 1, \ldots, d\,.
$$
Let $Z(u)$ and $Q_\alpha(u)$ be as in Theorem \ref{th:2}, and let
\begin{eqnarray*}
S_{\alpha, \delta}(u) & \defeq & 
\Expect{ \Big|Z(u)-Q_\alpha(u)\Big|^\delta\,\Big|\, Z(u) \le Q_\alpha(u) }\quad \mbox{and}
\\
S^\ast_{\alpha, \delta}(u) & \defeq & 
\Expect{ \Big|Z(u)-Q_\alpha(u)\Big|^\delta\,\Big|\, Z(u) \ge Q_\alpha(u) }\,,\ u \in U\,.
\end{eqnarray*}
If the density 
$t \mapsto |u_1|^{-1}\, \Expect{ \phi\left( u_1^{-1}\left( t -\ts{ \sum_{j=2}^d u_j\,X_j}\right),
X_2, \ldots, X_d\right) }$ of $Z(u)$ is positive at $t = Q_\alpha(u)$,
then for each $i=1,\ldots,d$ the partial derivatives of $S_{\alpha, \delta}$ and
$S^\ast_{\alpha, \delta}$ with respect
to $u_i$ exist and are continuous in some neigbourhood of $u$. They can be computed by
\begin{eqnarray}
\frac{\partial S_{\alpha, \delta}(u)}{\partial u_i} & = &
\delta\,  \Expect{\Big( \frac{\partial Q_\alpha(u)}{\partial u_i}-X_i\Big)
 \Big|Z(u)-Q_\alpha(u)\Big|^{\delta-1}\,\Big|\, Z(u) \le Q_\alpha(u)}\quad \mbox{and}
  \label{eq:13}\\
\frac{\partial S^\ast_{\alpha, \delta}(u)}{\partial u_i} & = &
\delta\,  \Expect{\Big(X_i - \frac{\partial Q_\alpha(u)}{\partial u_i}\Big)
 \Big|Z(u) - Q_\alpha(u)\Big|^{\delta-1}\,\Big|\, Z(u) \ge Q_\alpha(u)},\ u \in U\,,
\nonumber
\end{eqnarray}
where the formulas for $\frac{\partial Q_\alpha(u)}{\partial u_i}$ are given in Theorem \ref{th:2}.
  \label{th:3}
\end{theorem}
\textbf{Proof.}
For any event $A$ define the indicator function $I(A,\omega)=I(A)$ by
(\ref{eq:2.0}). Fix $\delta \ge 1$ and
let $G(z,u) \defeq \Expect{|Z(u)-z|^\delta\,I(Z(u) \le z)}$. 
Then the proof for e.g. (\ref{eq:13}) can be based on the representation
$$
G(z,u) \quad = \quad \delta \lint_{-\infty}^0 
\P{Z(u) \le z +t }\,(-t)^{\delta -1}\, d t\,,
$$
by using Lemma \ref{le:1} and Theorem A.(9.1) from  \cite{Durrett}. We omit the details.
\hfill$\Box$

Casually, one might be more interested in conditional moments of the underlying 
random variable itself than in moments of the difference between the random variable and a quantile.
The following corollary to Theorem \ref{th:3} covers this case. Note that
in contrast to our notation some people call solely the quantities $T_{\alpha, 1}(u)$
and $T^\ast_{\alpha, 1}(u)$,
defined in (\ref{eq:16}) and (\ref{eq:16a}) respectively, expected shortfall.

\begin{corollary}
Let an integer $n \ge 1$ be fixed. Let $(X_1, \ldots, X_d)$, $\phi$, $U$, and 
$\alpha$ be as in Theorem \ref{th:2} and assume 
additionally
$$
\Expect{\mid X_i\mid^n} < \infty,\quad i = 1, \ldots, d\,.
$$
Define $Z(u)$ and $Q_\alpha(u)$ as in Theorem \ref{th:2} and let
\begin{eqnarray}
T_{\alpha, n}(u) & \defeq & 
\Expect{ Z(u)^n\,\Big|\, Z(u) \le Q_\alpha(u) }\quad\mbox{and}
\label{eq:16}\\
T^\ast_{\alpha, n}(u) & \defeq & 
\Expect{ Z(u)^n\,\Big|\, Z(u) \ge Q_\alpha(u) }\,, \quad u \in U\,.
\label{eq:16a}
\end{eqnarray}
If the density 
$t \mapsto |u_1|^{-1}\, \Expect{ \phi\left( u_1^{-1}\left( t -\ts{ \sum_{j=2}^d u_j\,X_j}\right),
X_2, \ldots, X_d\right) }$ of $Z(u)$ is positive at $t = Q_\alpha(u)$,
then for each $i=1,\ldots,d$ the partial derivatives of $T_{\alpha, n}$ and 
$T^\ast_{\alpha, n}$ with respect
to $u_i$ exist and are continuous in some neighbourhood of $u$. They can be computed by
\begin{eqnarray}
\frac{\partial T_{\alpha, n}(u)}{\partial u_i} & = &
n\,  \Expect{X_i \,Z(u)^{n-1}\,\Big|\, Z(u) \le Q_\alpha(u)}\quad\mbox{and}
\nonumber\\
\frac{\partial T^\ast_{\alpha, n}(u)}{\partial u_i} & = &
n\,  \Expect{X_i \,Z(u)^{n-1}\,\Big|\, Z(u) \ge Q_\alpha(u)},\ u \in U\,.
  \label{eq:17a}
\end{eqnarray}
  \label{co:1}
\end{corollary}
\textbf{Proof.}
Follows 
through the representations 
\begin{eqnarray*}
  \Expect{ Z(u)^n \,\big|\, Z(u) \le Q_\alpha(u)} & = &
\sum_{k=0}^n {n \choose k} \, \Expect{ (Z(u)-Q_\alpha(u))^k \,\big|\, Z(u) \le Q_\alpha(u)}
\,Q_\alpha(u)^{n-k}\quad \mbox{and}\\
  \Expect{ Z(u)^n \,\big|\, Z(u) \ge Q_\alpha(u)} & = &
\sum_{k=0}^n {n \choose k} \, \Expect{ (Z(u)-Q_\alpha(u))^k \,\big|\, Z(u) \ge Q_\alpha(u)}
\,Q_\alpha(u)^{n-k}
\end{eqnarray*}
from Theorem \ref{th:3}. \hfill$\Box$

Equation (\ref{eq:17a}) in case $n=1$ was derived in \cite{Scaillet00} for random
variables $X_1, \ldots, X_d$ with a joint density.

\section{Expected shortfall as risk measure}
\label{sec:5}

Since quantiles (or the value-at-risk) as risk measures have some severe deficiencies
(cf. \cite{ADEH99} or \cite{Root_cklu}),
in the scientific literature other risk measures are preferred. One among those
risk measures is the Conditional Value-at-Risk (CVaR, expected shortfall with $n=1$ in
Corollary \ref{co:1}) because of its close relationship to the ``coherent''
risk measures introduced in \cite{ADEH99}.
 
Indeed, there seems to be some confusion in the literature whether 
in general CVaR is a coherent risk measure or not. For instance, in \cite{ADEH99} and
\cite{Delbaen} is argued that CVaR is not a coherent risk measure, whereas
\cite{Pflug00} says it is. The discrepancy between these statements is easy 
to explain since the definitions of CVaR the different authors used are not
identical. We will not examine here which definition is more useful, but
will show that even in the situation of the CVaR defined to be an 
elementary expectation -- as in   \cite{ADEH99} and \cite{Delbaen} --,
it enjoys to a great extent the coherence properties. We start with
an elementary but useful lemma.

\begin{lemma}
  \label{lem1.6}
Let $Y$ be a random variable in $L_1(\Omega, \mathfrak{F}, \mbox{P})$ and
$y\in\bbr$ such that $\P{ Y \le y }$ is positive. Then for all $F \in \mathfrak{F}$
with $\P{F} \ge \P{ Y\le y}$ we have
\begin{equation}
  \label{eq:1.6}
  \Expect{ Y\,|\,F}\quad \ge\quad \Expect{ Y\,|\,Y \le y}\,.
\end{equation}
\end{lemma}
\textbf{Proof.}
(\ref{eq:1.6}) is trivial in case $\P{ F \cap \{Y \le y\}} =
0$. Hence assume $\P{ F \cap \{Y \le y\}}$ to be positive. 
Define the indicator function $I(A)$ of the event $A \in \mathfrak{F}$ as in
(\ref{eq:2.0}).
We then obtain\\
\begin{minipage}[b]{13cm}
\begin{eqnarray*}
\Expect{ Y\,|\,Y \le y} & = & y\  +\ \frac{\Expect{(Y-y)\, I(\{Y\le y\}\cap F)} + 
                         \Expect{(Y-y)\, I(\{Y\le y\}\cap \Omega \bs F)} }{\P{Y\le y}} \\[0.5ex]
& \le & y\ +\  \Expect{Y - y\,|\,\{Y\le y\} \cap F}\,\P{F\,|\,Y\le y} \\[0.5ex]
& \le & y\ +\ \Expect{Y - y\,|\,\{Y\le y\} \cap F}\,\P{Y\le y\,|\,F} \\[0.5ex]
& \le & y\ +\ \frac{\Expect{(Y-y)\, I(\{Y\le y\}\cap F)} + 
                         \Expect{(Y-y)\, I(\{Y > y\}\cap  F)} }{\P{F}}\\[0.5ex]
& = & \Expect{Y\,|\,F}\,. 
\end{eqnarray*}
\end{minipage}
\hfill
\parbox[t]{1cm}{\hfill$\Box$}

Properties (i) to (iv) in the following proposition are just the constituting
properties of coherent risk measures (cf. \cite{Delbaen}). From this point
of view, Proposition \ref{pr:2} says that CVaR restricted to sets of
continuous random variables is in fact a coherent risk measure.

\begin{proposition}
Let $\alpha \in (0,1)$ be fixed. Assume that $\M$ is a convex cone 
(i.e. $X, Y \in \M, h > 0 \Rightarrow X+Y \in \M, h\, X \in \M$)
in $L_1(\Omega, \mathfrak{F}, \mbox{P})$.
Assume further that for each $X \in \M$ we have $\P{X \le Q_\alpha(X)} =  \alpha$ 
and $X + a \in \M$ for all $a \in \bbr$.

For $X \in \M$ define 
$$
\rho(X) \quad \defeq \quad  -\,\Expect{X\,|\, X \le Q_\alpha(X)}\,.
$$
Then we have
\begin{enumerate}
\item $\rho$ is monotonous, i.e. $X, Y \in \M, X \le Y \mbox{ a.s. } \Rightarrow 
\rho(X) \ge \rho(Y)$.
\item $\rho$ is subadditive, i.e. $X, Y \in \M \Rightarrow \rho(X+Y) \le \rho(X) + \rho(Y)$.
\item $\rho$ is positively homogeneous, i.e. $X \in \M, h > 0 \Rightarrow \rho(h\,X) = h\,\rho(X)$.
\item $\rho$ is translation invariant, i.e. $X \in \M, a \in \bbr 
\Rightarrow \rho(X + a) = \rho( X) - a$.
\end{enumerate}
  \label{pr:2}
\end{proposition}
\textbf{Proof.}
(iii), (iv) are trivial.
Concerning (i), by Lemma \ref{lem1.6} for $X \le Y$ we have
$$
  \rho(X) \ = \ -\,\Expect{X\,|\, X \le Q_\alpha(X)}
          \ \ge \ -\,\Expect{X\,|\, Y \le Q_\alpha(Y)}
          \ \ge \ -\,\Expect{Y\,|\, Y \ge Q_\alpha(Y)}
          \ = \ \rho(Y)\,.
$$
Similarly, concerning (ii):
\begin{eqnarray*}
  \rho(X + Y ) & = & -\,\Expect{X\,|\, X+Y \le Q_\alpha(X+Y)}\ -\ \Expect{Y\,|\, X+Y \le Q_\alpha(X+Y)}\\
               & \le & -\,\Expect{X\,|\, X \ge Q_\alpha(X)}\ -\ \Expect{Y\,|\, Y \ge Q_\alpha(Y)} \\
               & = & \rho(X)\ +\ \rho(Y)\,.
\end{eqnarray*}
This completes the proof.\hfill $\Box$

Note that the context specified by Lemma \ref{le:2.2} (and by Assumption \ref{as:2.3})
fits into the
assumptions for Proposition \ref{pr:2}. We state this fact formally
in the subsequent example.

\begin{example}
Suppose that $(X_1, \ldots, X_d)$ is an $\bbr^d$-valued random
vector with conditional density $\phi$ of $X_1$ given
$(X_2, \ldots, X_d)$. Assume that $\Expect{|X_i|}< \infty, i = 1, \ldots, d$. 
Let $U \subset (0,\infty) \times \bbr^{d-1}$ be an open convex cone and
define
$$
\M \quad \defeq \quad \{ a + \sum_{i=1}^d u_i\,X_i\,|\, a \in \bbr, u \in U \}\,.
$$
Then $\M$ satisfies the conditions of Proposition \ref{pr:2}. 
\end{example}

\bibliographystyle{plain}
\bibliography{quant}

\end{document}